\documentclass[12pt]{article}
\usepackage[centertags]{amsmath}
\usepackage{amsfonts}
\usepackage{amssymb}
\usepackage{amsthm}
\usepackage{newlfont}
\usepackage{graphics}
\usepackage{mathrsfs}
\usepackage{hyperref}

\newcommand{\vm}[1]{\mbox{\boldmath$#1$}}
\newcommand{\Prob}[1]{\Pr \left\{ #1 \right \}}
\newcommand{\eps}{\varepsilon}
\newcommand{\sqep}{\sqrt{\varepsilon}}

\newcommand{\ds}{\displaystyle}
\newcommand{\beq}{\begin{eqnarray}}
\newcommand{\beqq}{\begin{eqnarray*}}
\newcommand{\eeq}{\end{eqnarray}}
\newcommand{\eeqq}{\end{eqnarray*}}
\newcommand{\x}{\mbox{\boldmath$x$}}
\newcommand{\y}{\mbox{\boldmath$y$}}
\newcommand{\z}{\mbox{\boldmath$z$}}
\newcommand{\w}{\mbox{\boldmath$w$}}

\newcommand{\e}{\mbox{\boldmath$e$}}

\begin{document}
\hsize=6 in  \textheight=9 in
\bibliographystyle{plain}
\title{ \textbf{THE EXIT PROBLEM IN OPTIMAL NON-CAUSAL ESTIMATION}}
\author {D. Ezri\thanks{Department of Electrical Engineering--Systems, Tel-Aviv University, Ramat-Aviv, Tel-Aviv 69978,
Israel. email: ezri@eng.tau.ac.il }   \and B. Z.
Bobrovsky\thanks{Department of Electrical Engineering--Systems,
Tel-Aviv University, Ramat-Aviv, Tel-Aviv 69978, Israel. email:
bobrov@eng.tau.ac.il} \and Z. Schuss\thanks{Department of
Mathematics, Tel-Aviv University, Ramat-Aviv, Tel-Aviv 69978,
Israel. email: schuss@post.tau.ac.il } }

\maketitle

\begin{abstract}
We study the phenomenon of loss of lock in the optimal non-causal
phase estimation problem, a benchmark problem in nonlinear
estimation. Our method is based on the computation of the asymptotic
distribution of the optimal estimation error in case the number of
trajectories in the optimization problem is finite. The computation
is based directly on the minimum noise energy optimality criterion
rather than on state equations of the error, as is the usual case in
the literature. The results include an asymptotic computation of the
mean time to lose lock (MTLL) in the optimal smoother. We show that
the MTLL in the first and second order smoothers is significantly
longer than that in the causal extended Kalman filter.

\end{abstract}


\section{Introduction}
In many applications in communication practice a random signal
$\x(t)$ is received through a noisy channel. The random signal $
 \x(t) \in \mathbb{R}^N$ is assumed to be a stochastic
process defined by an It\^o  stochastic  differential equation (SDE)
\cite{Liptser}
\begin{equation}\label{model1}
    d\x=\vm{m}(\x,t)\,dt +\vm{\sigma}(\x,t)\, d\w,
\end{equation}
where $\w(t)$ is a  vector of standard Brownian motions. The
measurements process
  $\y(t) \in \mathbb{R}^M$, which is the output of the noisy channel, is modeled by another
  SDE
\begin{equation}\label{measure1}
    d\y=\vm{g}(\x,t)\,dt+\sqrt{N_0/2}\, d\vm{v},
\end{equation}
where $N_0$ measures the channel noise intensity and $\vm{v}(t)$
is another vector of Brownian motions, independent of $\w(t)$.  We
further assume that the functions $\vm{m}(\cdot,\cdot),\,
\vm{\sigma}(\cdot,\cdot)$ and $\vm{g}(\cdot,\cdot)$ satisfy
standard regularity conditions such that
\eqref{model1},\eqref{measure1} possess a strong unique solution.

When the optimality criterion is minimum square error, the optimal
filtering problem is to construct the causal estimator
$\hat{\x}(t)=E\left[\x(t)\,|\,\y(s)\right ]$ of $\x(t)$, where
$0\leq s \leq t$
 \cite{VanTrees}. The optimal fixed interval smoothing
problem is to construct the non-causal estimator
$\hat{\x}(t)=E\left[\x(t)\,|\,\y(s), 0\leq s \leq T\right]$, where
$T$ is the length of the interval, and $t<T$. The optimal fixed lag
smoothing problem is  to construct the non-causal estimator
$\hat{\x}(t)=E\left[\x(t)\,|\,\y(s), 0\leq s \leq t+\Delta\right]$,
where $T$ is the length of the interval, and $t+\Delta<T$. In many
applications delay in the estimation of the signal is not
permissible, as for example in closed-loop control, radar tracking
systems, and so on. There are, however, interesting cases, where
certain delay is permissible, as for example in communication
systems, as extensively practiced in coding
\cite{ViterbiCoding,Proakis}.

Smoothers are used because their performance is superior to all
causal filters, with respect to the same optimality criterion
\cite{meditch}. In linear estimation theory, where the optimality
criterion is minimum mean square error, the error variance of the
optimal smoother is smaller than that of the optimal filter
\cite{kailath22}.

Optimal estimators are usually infinite-dimensional and therefore
have no finite-dimensional realizations, so that suboptimal
estimators have been proposed to approximate the optimal ones by a
system of SDEs, driven by the measurements
\cite{kushner,Rauch1,Rauch2}. The phase-locked-loop (PLL), which is
a realization of the extended Kalman filter (EKF) \cite{Snyder}, is
a nonlinear suboptimal causal estimator of the carrier phase in
various communications systems \cite{Lindsey}. A well known effect
in such PLL demodulators is the cycle slip phenomenon that consists
in occasional sudden changes of size $2\pi\,n \ (n=\pm 1, \pm
2,\dots)$ in the phase estimation error \cite{Benzi82}. Obviously,
the mean time between cycle slips, known as the mean time to lose
lock (MTLL), decreases with the noise intensity and causes sharp
degradation in the performance of the filter and to the formation of
a performance threshold ~\cite{Viterbi,Schuss_book,Schuss_article}.

Considerable effort was put into the computation of the MTLL in
causal estimators \cite{Lindsey,Ryter,Tausworth,Viterbi},
including the singular perturbation method
\cite{Schuss_article,Schuss_book} and large deviations theory
\cite{WentzellFreidlin,Zeitouni}. However, the phenomenon of loss
of lock  in smoothers has never been addressed, despite the
extensive study of linear and nonlinear smoothers  in the
literature
\cite{Rauch1,Rauch2,Rauch3,Lee,Weaver,bryson,BellmanKalaba,Sage,SageEwing,kailath22,Leondes,Anderson}.
The objective of this paper is to provide the missing theory,
estimate the MTLL in the optimal smoother, and compare it with
that in the casual PLL. Specifically, we compute the asymptotic
distribution of the optimal estimation error in case the number of
trajectories in the optimization problem is finite. We identify
the contribution of error trajectories to the minimum noise energy
(MNE) cost functional and recast the problem in terms of order
statistics. The asymptotic expression for the MTLL in the smoother
is similar to that resulting from the Wentzell-Freidlin theorem
for causal systems, with a new functional. Applying our method to
standard phase models, we show that the MTLL in the optimal
smoother is significantly longer than that in the PLL.

\section{The mathematical model}

The general equations of a scaled phase tracking system consist of
the linear model of the phase $\x(t)=[x(t), x_2(t),...,x_N(t)]^T$
\cite{Schuss_book}
 \beq\label{GenralScaledSystemx}
 \dot{\x}&=&\vm{A}\x+\sqep\,\vm{B}
\,\dot{\w}
 \eeq
and the nonlinear model of the noisy measurements
$\y(t)=[y_s(t),y_c(t)]^T$
\begin{eqnarray}\label{GenralScaledSystemy}
 \y=\vm{h}(\x) + \sqep\,\dot{\vm{v}},
\end{eqnarray}
with
\begin{equation*}
\vm{h}(\x) =\left(\begin{array}{c}
    \sin x \\
    \cos x
  \end{array}   \right).
\end{equation*}
The dimensionless parameter $\eps$ is assumed small in the case of a
low noise channel \cite{Schuss_article}.

A fixed-interval minimum noise energy (MNE) estimator
$\hat{\x}(\cdot)$ for $\x(\cdot)$  is the minimizer of the cost
functional \cite{bryson}
\begin{eqnarray}\label{defineJ1}
 J[\z(\cdot)] & \equiv &
\int_{0}^{T} \left[ \left | \y-\vm{h}(\z) \right | ^2
+|\vm{\zeta}|^2\right]\,dt,
\end{eqnarray}
with the equality constraint
\begin{eqnarray}\label{DefineEqCon}
\dot{\z}=\vm{Az}+\vm{B\zeta},
\end{eqnarray}
that is,
\begin{equation}
\label{DefineXHat}
 \hat{\x}(\cdot)\equiv \textrm{arg}\min_{\z(\cdot)}J[\z(\cdot)].
\end{equation}

Note that the integral $J[\x(\cdot)]$ contains the white noises
$\dot{\w}(t),\ \dot{\vm{v}}(t)$, which are not square integrable.
To remedy this problem, we begin with a model in which the white
noises $\dot{\w}(t),\ \dot{\vm{v}}(t)$ are replaced with square
integrable wide band noises, and at the appropriate stage of the
analysis, we take the white noise limit (see below).

In contrast to nonlinear filtering, where the locked state is a
local attractor for the error dynamics \cite{Schuss_book}, and
escaping it corresponds to loss of lock, there is no dynamics, and
therefore no attractors for smoothers. Thus, we have to define
cycle slip events in a different manner than hitting the boundary
of the domain of attraction. We define the estimation error
$\e(t)=[e(t), e_2(t), \ldots, e_N(t)]^T$ as
\begin{equation}\label{SDefinition}
    \e(t)=\hat{\x}(t)-\x(t),
\end{equation}
and we say that a cycle-slip has occurred in the time interval
$(t_0, t_0+\Delta t), \ \Delta t << 1$, if and only if the
estimation error $\e(t)$ vanishes at at some $t_0-T_1$, reaches the
point $[2\pi n,\underline{\vm{0}}]^T, \ (n=\pm1, \pm 2,...)$, at
some later time $t_0+T_2$, and $e(\tilde{t})=\pi n$, where
$\tilde{t} \in (t_0,t_0+\Delta t)$. The time $T_s=T_1+T_2$ is the
slip duration, satisfying $T_s<<T$. We define in the space of
continuous functions $\mathscr{C}^N_{[0,T]}$ the set
$\mathscr{C}^N(t_0)$ of all continuous trajectories $\e(\cdot)$ that
slip in the interval $(t_0, t_0+\Delta t)$. Thus
\begin{equation}\label{CSdef}
     \Prob{\mbox{slip in} \ (t_0, t_0 +\Delta t)}=\Prob{\e(\cdot)\in
     \mathscr{C}^N(t_0)}.
\end{equation}

\begin{figure}
\centering
\resizebox{!}{9cm}{\includegraphics{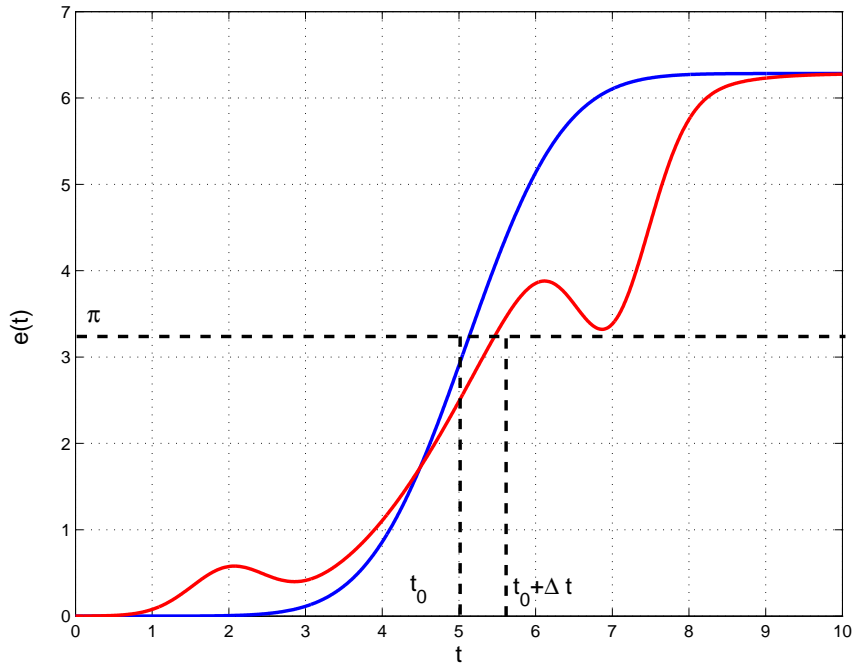}}
\caption{\label{FunctionsInC1}An example of two trajectories in
$\mathcal{C}^1(t_0=5)$}
\end{figure}

For small values of $\eps$ the cycle-slip event is a rare large
deviation from the original trajectory, and therefore $\hat{\x}(t)$
is in the vicinity of $\x(t)$ before the cycle-slip, and in the
vicinity of $\hat{ \x}(t)+[2 \pi n,\underline{\vm{0}}]^T$ after the
slip. Thus, we can define the beginning and the end of the
cycle-slip by the instants when $\e(t)$ reaches the origin and
$[2\pi n,\underline{\vm{0}}]^T$, respectively. An example of two
trajectories in $\mathscr{C}^1(t_0=5)$ is given in Figure
\ref{FunctionsInC1}.

\section{The MTLL in the optimal smoother}
The Wentzell-Freidlin theorem \cite{Zeitouni,WentzellFreidlin} and
the singular perturbation method \cite{Schuss_book} for asymptotic
evaluation of the MTLL are concerned with stochastic processes
satisfying a stochastic differential equation with a unique
solution. In contrast, the dynamics of the optimal smoother,
derived from the EL equations, form a two-point boundary-value
problem which has no unique solution. Therefore the
Wentzell-Freidlin and the singular perturbation method seem
inappropriate for the computation of the MTLL in a smoother. It
appears that this computation calls for a different approach.

The first step toward an asymptotic calculation of the MTLL in
smoothers is the computation of the asymptotic distribution of the
estimation error in case the number of trajectories in the
optimization problem is finite. We investigate the cost functional
of deterministic error trajectories that deviate from the original
trajectory $\x(t)$. We augment $\x(t)$ with the set of the $N_T$
trajectories $\vm{r}_i(t)=[r_i(t), r_i^{[2]}(t),
\ldots,r_i^{[N]}(t)] \in \mathscr{C}^N_{[0,T]},\ i \in
[1,\ldots,N_T]$. The trajectories $\x(t)+\vm{r}_i(t)$ are
admissible in the optimization problem \eqref{defineJ1},
\eqref{DefineEqCon}, only if the trajectories $\vm{r}_i(t)$
satisfy
\begin{equation}\label{defineR}
\dot{\vm{r}}_i=\vm{A}\vm{r}_i+ \vm{B}\vm{u}_i.
\end{equation}

We define the difference
\begin{equation}\label{definedeltaJ}
\Delta J[\x(\cdot),\vm{r}_i(\cdot)]\stackrel\triangle
=J[\x(\cdot)+\vm{r}_i(\cdot)]-J[\x(\cdot)]
\end{equation}
and substitute \eqref{defineJ1}, \eqref{DefineEqCon} and
\eqref{defineR} in \eqref{definedeltaJ} to obtain

\begin{eqnarray}\label{Delta}
\begin{array}{lll}
 \Delta J[\x(\cdot),\vm{r}_i(\cdot)]&=&
\ds\int_{0}^{T}{ \left[ 4 \sin^2 \left ( \ds\frac{r_i}{2}\right )+
|\vm{u}_i|^2 \right ] }\,dt \\&&\\
 & + & \sqrt{4\,\eps} \ds\int_{0}^{T}{ \left[ \sin x -
\sin(x+r_i) \right] }\,dv_1(t)  \\&&\\
 & + & \sqrt{4\,\eps}
\ds\int_{0}^{T}{ \left[ \cos x
 -\cos(x+r_i) \right] }\,dv_2(t) \\&&\\
& + & \sqrt{4\,\eps} \ds\int_{0}^{T} \vm{u}_i^T \,d\w(t).
\end{array}
\end{eqnarray}
At this point, we take the white noise limit in the wide band noises
so the stochastic integrals in \eqref{Delta} become It\^o integrals.
Collecting them into a single It\^o integral leads to
\begin{eqnarray}
\begin{array}{lll}
\Delta J[\x(\cdot),\vm{r}_i(\cdot)]& =& \ds\int_{0}^{T}\left [
4\sin^2\left(\ds\frac{r_i}2\right)+|\vm{u}_i|^2\right ]\,dt \\ &&
\\
&+& \sqrt{4\,\eps} \ds\int_{0}^{T}
\sqrt{4\sin^2\left(\ds\frac{r_i}2\right)+|\vm{u}_i|^2}\,\,d\tilde{v}_i(t),
\end{array}
\end{eqnarray}
where $\tilde{v}_i(t)$ is a standard Brownian motion that depends on
$\vm{v}(t)$ and $\w(t)$.

Note that although the values of the random variable $\Delta
J[\x(\cdot),\vm{r}_i(\cdot)]$ depend on the trajectories of $\x(t)$
and $\y(t)$ through $\vm{v}(t)$ and $\w(t)$, the probability law of
$\Delta J[\x(\cdot),\vm{r}_i(\cdot)]$ depends only on the
trajectories  $\vm{r}_i(t)$. Therefore, we abbreviate notation to
$\Delta J[\vm{r}_i(\cdot)]$. We note further that
 $\Delta J[\vm{r}_i(\cdot)]$ is  a Gaussian random variable
with expectation
\begin{equation}
\mbox{E}\Delta J[\vm{r}_i(\cdot)]=m_i,
\end{equation}
where
\begin{equation}\label{Psidef}
m_i \stackrel{\triangle}{=} \int_{0}^{T}\left [
4\sin^2\left(\ds\frac{r_i}2\right)+|\vm{u}_i|^2\right ]\,dt,
\end{equation}
and  variance
\begin{equation}
\mbox{Var}\{\Delta J[\vm{r}_i(\cdot)]\} = 4\eps \,m_i.
\end{equation}
Furthermore, we compute the covariance
\begin{eqnarray}
\begin{array}{lll}
\sigma_{i\,j} & =& \mbox{E}\left ( \Delta J[\vm{r}_i(\cdot)]-m_i
\right )\left ( \Delta J[\vm{r}_j(\cdot)]-m_j \right )\\ && \\
&=& 4 \eps \ds \int_0^T \vm{u}_i^T \vm{u}_j\, dt + 4 \eps \ds
\int_0^T \left [ 1+ \cos(r_i-r_j)-\cos r_i -\cos r_j    \right ]\,
dt.
\end{array}
\end{eqnarray}
Note  that if the supports of $\vm{r}_i(\cdot)$ and
$\vm{r}_j(\cdot)$ are disjoint, the cost functionals $\Delta
J[\vm{r}_i(\cdot)], \,$ $\Delta J[\vm{r}_j(\cdot)]$ are not
correlated.  We conclude that the random variables $\Delta
J[\vm{r}_i(\cdot)],\ i=1\ldots N_T$, form a Gaussian random vector
with distribution
\begin{eqnarray}
\left [ \begin{array}{c}
  \Delta J[\vm{r}_1(\cdot)] \\ \\
  \Delta J[\vm{r}_2(\cdot)] \\ \\
  \vdots \\ \\
  \Delta J[\vm{r}_{N_T}(\cdot)] \\
\end{array}
\right ] \sim N \left \{ \left [
\begin{array}{c}
  m_1 \\ \\
  m_2 \\ \\
  \vdots \\ \\
  m_{N_T} \\
\end{array}
\right ] \ ; \ 4\,\eps \left [ \begin{array}{cccc}
  m_1 & \sigma_{1\,2} & \ldots & \tilde{\sigma}_{1\,N_T} \\ &&& \\
    \tilde{\sigma}_{1\,2} & m_2 &  &  \\&&&\\
  \vdots &  & \ddots &  \\&&&\\
   \tilde{\sigma}_{1\,N_T}&  &  & m_{N_T} \\
\end{array}\right ]
\right \},
\end{eqnarray}
where $\tilde{\sigma}_{i j}=\ds\frac{1}{4\,\eps}\sigma_{i j}$.

When considering a finite number of error trajectories
$\vm{r}_i(\cdot),\, i=1,\ldots, N_T$ in the optimization problem.
The estimator error trajectory $\e_{N_T}(\cdot)$ minimizes the
cost functional $\Delta J[\vm{r}_j(\cdot)]$,
\begin{eqnarray}\label{ErrorProbToA}
\Prob{\e_{N_T}(\cdot)=\vm{r}_k(\cdot)}=\Prob{\Delta
J[\vm{r}_k(\cdot)]<\Delta J[\vm{r}_j(\cdot)] \,\mbox{for all} \ j
\neq k }.
\end{eqnarray}
Thus the problem of minimization has been recast in the language
of order statistics. The probability on the right side of
\eqref{ErrorProbToA} is difficult to calculate, so we pursue the
distribution of the error trajectory $\e_{N_T}(\cdot)$ in the
limit of small $\eps$. We assume that for each $k$ there exists an
interval $A_k$ such that
\begin{eqnarray}\label{BasicInequality}
\mbox{E}\left \{ \Delta J[\vm{r}_j(\cdot)]\,|\, \Delta
J[\vm{r}_k(\cdot) \in A_k ] \right \} > \Delta J[\vm{r}_k(\cdot)]
\quad\mbox{for all}\quad  j \neq k.
\end{eqnarray}
Cram\'{e}r's theorem for Gaussian vectors \cite{Zeitouni,Bucklew}
implies that in the limit of small $\eps$
\begin{eqnarray}\label{AsymptoticA}
\begin{array}{lll}
&&\ds \lim_{\eps \to 0} \eps \log_e \Prob{\Delta J [\vm{r}_k(\cdot)]
<\Delta J [\vm{r}_j(\cdot)] \quad\mbox{for all}\quad  j \neq k
}=\\&&\\
&&\ds\lim_{\eps \to 0} \eps \log_e \Prob{\Delta J [\vm{r}_k(\cdot)]
\in A_k }.
\end{array}
\end{eqnarray}
Next, we evaluate the conditional expectation
\begin{eqnarray}
\mbox{E}\left \{ \Delta J[\vm{r}_j(\cdot)]\,|\, \Delta
J[\vm{r}_k(\cdot)] \right \}.
\end{eqnarray}
Since $\Delta J[\vm{r}_j(\cdot)],\,\Delta J[\vm{r}_k(\cdot)] $ are
jointly Gaussian,
\begin{eqnarray}
\begin{array}{lll}
\mbox{E}\left \{ \Delta J[\vm{r}_j(\cdot)]\,|\, \Delta
J[\vm{r}_k(\cdot)] \right \}&=& m_j+\ds\frac{\sigma_{j k}}{4\,\eps
m_k} \left ( \Delta J[\vm{r}_k(\cdot)]-m_k   \right ) \\ && \\
& = & m_j+\ds\frac{\tilde{\sigma}_{j k}}{m_k} \left ( \Delta
J[\vm{r}_k(\cdot)]-m_k \right ).
\end{array}
\end{eqnarray}
In order to determine the interval $A_k$, defined in
\eqref{BasicInequality}, we derive the set of $N_T-1$ inequalities
\begin{eqnarray}\label{InequalityBasic}
\mbox{E}\left \{ \Delta J[\vm{r}_j(\cdot)]\,|\, \Delta
J[\vm{r}_k(\cdot) ] \right \} > \Delta J[\vm{r}_k(\cdot)] \
\mbox{for all} \ j \neq k.
\end{eqnarray}
The interval $A_k$, which is the range of values of $\Delta J
[\vm{r}_k(\cdot)]$ satisfying \eqref{InequalityBasic} is defined by
\begin{eqnarray}\label{InequalForA}
\mathop{\max}\limits_{\vm{r}_j(\cdot) \in A^-}
\frac{m_j-\tilde{\sigma}_{j k} }{ 1 - \ds\frac{\tilde{\sigma}_{j
k}}{m_k}}<\Delta J
[\vm{r}_k(\cdot)]<\mathop{\min}\limits_{\vm{r}_j(\cdot) \in A^+ }
\frac{m_j-\tilde{\sigma}_{j k} }{ 1 - \ds\frac{\tilde{\sigma}_{j
k}}{m_k}},
\end{eqnarray}
where $A^+$ is the set of all trajectories $\vm{r}_j(\cdot)$ such
that $m_j<m_k,\ \tilde{\sigma}_{j k}>m_j$, and $A^-$ is the set of
all trajectories $\vm{r}_j(\cdot)$ such that $m_j>m_k,\
\tilde{\sigma}_{j k}>m_k$. Note that the supremum and infimum,
over all continuous trajectories, of the leftmost and rightmost
sides of \eqref{InequalForA}, respectively, is $-m_k$. This means
that as $N_T$ increases and the trajectories
$\x(\cdot)+\vm{r}_i(\cdot)$ are sampled from ${\mathscr
C}_{[0,T]}$ according to their {\em a priori} distribution
\eqref{GenralScaledSystemx}, the interval $A_k$ narrows.
Specifically, for any $\delta,\, \tilde{\delta}>0$ there is a
sufficiently large $N_T$ such that
$\Prob{A_k\not\subset(-m_k-\delta,-m_k+\delta)}<\tilde{\delta}$.
Combining \eqref{ErrorProbToA} and \eqref{AsymptoticA}, we
conclude that for all small $\eps>0$ and every sufficiently small
$\delta$, such that $0<\delta<\eps$, there is a sufficiently large
$N_T$ such that
\begin{eqnarray}\label{DistribFinite}
\begin{array}{lll}
\Prob{\e_{N_T}(\cdot) =
\vm{r}_k(\cdot)}&\asymp&\Prob{-m_k-\delta<\Delta
J[\vm{r}_k(\cdot)]<-m_k+\delta} \\ && \\
&\asymp& 2\delta\exp \left \{ -\ds \frac{4\,m_k^2}{8 \eps \,
m_k} \right \} \\ && \\
&\asymp& 2\delta\exp \left\{  -\ds \frac{1}{2 \eps} \ds \int_0^T
\left [ 4 \sin^2 \left (\ds \frac{r_k}{2}\right ) + |\vm{u}_k|^2
\right] \, dt \right\}.
\end{array}
\end{eqnarray}

Based on the distribution of the estimation error
$\vm{e}_{N_T}(\cdot)$ in the case of a finite number of error
trajectories $N_T$ \eqref{DistribFinite}, the probability that
$\vm{e}_{N_T}(\cdot)$ is in any set $\mathscr{A}$ in
$\mathscr{C}^N_{[0,T]}$ is
\begin{eqnarray}
\begin{array} {lll}
\Prob{\vm{e}_{N_T}(\cdot) \in
\mathscr{A}}&=&\ds\sum_{\vm{r}_k(\cdot) \in \mathscr{A}}
\Prob{\vm{e}_{N_T}(\cdot)=\vm{r}_k(\cdot)} \\ &&
\\
& \asymp &\ds\sum_{\vm{r}_k(\cdot) \in \mathscr{A}} 2\delta\exp
\left\{  -\ds \frac{1}{2 \eps} \ds \int_0^T \left [ 4 \sin^2 \left
(\ds \frac{r_k}{2}\right ) + |\vm{u}_k|^2 \right] \, dt \right\}.
\end{array}
\end{eqnarray}
Applying Laplace's method for sums of exponentials with large
parameter $\ds \frac{1}{\eps}$ \cite{Bender}, we obtain the
asymptotic expression
\begin{eqnarray}
\Prob{\vm{e}_{N_T}(\cdot) \in \mathscr{A}} \asymp  \exp \left\{-
\ds \frac{1}{2 \eps} \mathop{\min} \limits_{\vm{r}_k(\cdot) \in
\mathscr{A}} \ds \int_0^T \left [ 4 \sin^2 \left (\ds
\frac{r_k}{2}\right ) + |\vm{u}_k|^2 \right] \, dt \right\}.
\end{eqnarray}

In the limit $\delta\to0$ (and $N_T \to \infty$), the probability
that the optimal estimation error $\vm{e}(\cdot)$ is in
$\mathscr{A}$ is found as
\begin{eqnarray}\label{AsympProbInA}
\Prob{\vm{e}(\cdot) \in \mathscr{A}} \asymp  \exp \left\{ -\ds
\frac{1}{2 \eps} \mathop{\inf} \limits_{\vm{r}(\cdot) \in
\mathscr{A}} \ds \int_0^T \left [ 4 \sin^2 \left (\ds
\frac{r}{2}\right ) + |\vm{u}|^2 \right] \, dt \right\},
\end{eqnarray}
subject to the equality constraint
\begin{eqnarray}\label{EqualityConstraintForR}
\dot{\vm{r}}=\vm{Ar}+\vm{B u},
\end{eqnarray}
where the infimum in \eqref{AsympProbInA} is taken over all
continuous trajectories $\vm{r}(\cdot) \in \mathscr{A}$. Note that
the small $\varepsilon$ approximation is taken before the limit
$\delta\to0.$

We turn now to the computation of the MTLL. First, we consider
time intervals $[0,T]$ much longer than the time constant of the
system, so that most of the time the system is in steady state. It
follows that the probability $\Prob{\mbox{slip in} \
(t_0,t_0+\Delta t)}$ is independent of $t_0$ for $t_0$ outside
intervals of fixed length (the time constant of the system) at the
endpoints $0$ and $T$. Therefore, in view of the regularity of the
pdf of the solution as a function of $t$ and the independence of
cycle slips in disjoint intervals (see above), for such $t_0$
\begin{eqnarray}\label{timelinear}
\begin{array}{l}
\Prob{\mbox{slip in} \ (t_0,t_0+2\Delta t)}=\\ \\
\Prob{\mbox{slip in} \ (t_0,t_0+\Delta t)}+\Prob{\mbox{slip in} \
(t_0+\Delta t,t_0+2\Delta t)}+o(\Delta t)= \\  \\
2\Prob{\mbox{slip in} \ (t_0,t_0+\Delta t)}+o(\Delta t).
\end{array}
\end{eqnarray}
Thus, $\Prob{\mbox{slip in} \ (t_0,t_0+\Delta t)}$ is nearly
linear in $\Delta t$.

Next, we note that for fixed $\Delta t$ \eqref{AsympProbInA}
implies that the slip probability satisfies
\begin{eqnarray}\label{SlipProbAsIntegral}
\begin{array}{lll}
 && \Prob{\mbox{slip in $(t_0,t_0+\Delta t)$}}\\ && \\
 &=&\Prob{\e(\cdot)\in
\mathscr{C}^N(t_0)} \\ && \\
 & \asymp & \exp \left\{ -\ds
\frac{1}{2 \eps} \mathop{\inf} \limits_{\vm{r}(\cdot) \in
\mathscr{C}^N(t_0)} \ds \int_0^T \left [ 4 \sin^2 \left (\ds
\frac{r}{2}\right ) + |\vm{u}|^2 \right] \, dt \right\},
\end{array}
\end{eqnarray}
subject to the equality constraint \eqref{EqualityConstraintForR}.
Equations (\ref{timelinear}) and  (\ref{SlipProbAsIntegral}) can
be written together as
\begin{eqnarray}\label{SlipProb}
\begin{array}{lll}
&& \\
& &\Prob{\mbox{slip in $(t_0,t_0+\Delta t)$}}=\\ && \\
&&(\Delta t+o(\Delta t)) \Omega(\varepsilon)\exp \left
\{-\mathop{\inf} \limits_{\vm{r}(\cdot) \in \mathscr{C}^N(t_0)}
\ds\frac{1}{2\eps}\ds\int_0^{T}\left[4\sin^2 \left( \ds\frac{r}{2}
\right)+|\vm{u}|^2 \right]\,dt\right \},
\\&&
\end{array}
\end{eqnarray}
subject to the equality constraint \eqref{EqualityConstraintForR},
where $\varepsilon\log\Omega(\varepsilon)\to0$ as
$\varepsilon\to0$.

Equipped with the slip probability \eqref{SlipProb}, we turn to the
evaluation of the MTLL in the optimal smoother. We introduce a {\em
renewal} (counting) process $\left\{ \mathscr{N}(t),\,t\geq
0\right\} $, a nonnegative integer-valued stochastic process that
counts the number successive cycle-slip events in the time interval
$(0,t]$ \cite{Karlin}. We assume that the time durations between
consecutive slips are  positive, independent, identically
distributed random variables. Based on the above assumptions, we
adopt the renewal formula \cite{Karlin}
\begin{eqnarray}\label{RenewalArgument}
\tau_{nc}=\frac{t}{\mbox{E} \mathscr{N}(t)},
\end{eqnarray}
where $\tau_{nc}$ is the MTLL in the non-causal estimator. In the
limit of $t \to \infty$, equation \eqref{RenewalArgument} gives
\begin{eqnarray}\label{TauAndNdot}
\tau_{nc}=\lim_{t \to \infty}\frac{t}{\mbox{E}
\mathscr{N}(t)}=\lim_{t \to
\infty}\frac{1}{\mbox{E}\dot{\mathscr{N}}(t)}.
\end{eqnarray}
Using the slip probability  \eqref{SlipProb} we obtain
\begin{eqnarray}\label{N_dot}
\begin{array}{lll}
\mbox{E}\dot{\mathscr{N}}(t)&=&\ds \lim_{\Delta t \to 0}
\mbox{E}\left \{
 \ds\frac{\mathscr{N}(t+\Delta t)-\mathscr{N}(t)}{\Delta t} \right \}
 \\ && \\
 &=& \ds\lim_{\Delta t \to 0}\ds\frac{\Prob{\mbox{slip in}\ (t_0,t_0+\Delta t)}}{\Delta
 t}\\&&\\
 &\asymp& \exp \left \{-\mathop{\inf} \limits_{\vm{r}(\cdot) \in
\mathscr{C}^N(t_0)} \ds\frac{1}{2\eps}\ds\int_0^{T}\left[4\sin^2
\left( \ds\frac{r}{2} \right)+|\vm{u}|^2 \right]\,dt\right \}.
\end{array}
\end{eqnarray}
Substituting \eqref{N_dot} in \eqref{TauAndNdot}, we obtain the
expression for the asymptotic MTLL, $\tau_{nc}$, in the optimal
MNE estimator
\begin{equation}\label{generalMTLL}
 \lim_{\eps \rightarrow 0}\eps \log_e \tau_{nc}=
 \mathop{\inf} \limits_{\e(\cdot) \in \mathscr{C}^N(t_0)}
\frac{1}{2}\int_0^{T}\left [ 4\sin^2 \left( \frac{e}{2}
\right)+|\vm{\xi}|^2 \right ] \,dt,
\end{equation}
subject to the equality constraint
\begin{eqnarray}\label{EqualityConstraintForS}
\dot{\e}=\vm{A}\e+\vm{B\xi}.
\end{eqnarray}

\section{The MTLL in the smoother with standard phase models}

We begin with the first order phase tracking system suggested by
\cite{Scharf} and \cite{Macchi}, in which the phase $x(t)$ is
modelled as a standard Brownian motion
\begin{eqnarray}\label{sys_first}
\begin{array}{lll}
&& \\
\dot{x} &=& \dot{w} \\
&& \\
y&=&\vm{h}(x)+\rho \, \dot{\vm{v}}, \\
&&
\end{array}
\end{eqnarray}
where $x(t),w(t)$ take values in $ \mathbb{R}^1$. The system
\eqref{sys_first} is scaled to the form of
\eqref{GenralScaledSystemx}, \eqref{GenralScaledSystemy} with $
\vm{A}=0,\vm{B}=1$ and $\eps=\rho$. A similar procedure is presented
in \cite{Benzi78}. Note that the CNR equals $\rho^{-2}/2$ in the
system \eqref{sys_first}.

In the first order case, the asymptotic expression for the MTLL in
the smoother \eqref{generalMTLL} becomes
\begin{equation}\label{MTLLFirst}
 \lim_{\eps \rightarrow 0}\eps \log_e \tau_{nc}=
 \mathop{\inf} \limits_{e(\cdot) \in \mathscr{C}^1(t_0)}
\frac{1}{2}\int_0^{T}\left [ 4\sin^2 \left( \frac{e}{2}
\right)+\dot{e}^2 \right ] \,dt,
\end{equation}
The variational problem on the right side of \eqref{MTLLFirst} is
solved analytically using the Hamilton-Jacobi-Belman equation
\cite{Kirk}. This leads to the asymptotic limit of the MTLL in the
first order non-causal estimator as
\begin{equation}\label{ourMTLL}
\lim_{\eps\rightarrow 0}\eps \log_e \tau_{nc}= \lim_{\rho
\rightarrow 0}\rho \log_e \tau_{nc}=8.
\end{equation}

In order to compare the MTLL in the optimal smoother \eqref{ourMTLL}
with that in the suboptimal PLL we construct the steady-state EKF
corresponding to the model \eqref{sys_first} \cite{Snyder}.
\begin{eqnarray}
\dot{\hat{x}}=\frac{\sigma}{\rho}\left( y_s \cos \hat{x} -y_c \sin
\hat{x} \right).
\end{eqnarray}
The differential equation of the causal EKF estimation error
$e(t)=\hat{x}(t)-x(t)$ is scaled to the form \cite{Benzi78}
\begin{equation}\label{Kalman_error}
    \dot{e}=-\sin e  +\sqrt{2\eps}\,\dot{v}.
\end{equation}
The asymptotic MTLL, $\tau_c$, in this simple analytical potential
case is known to be \cite{Schuss_book}
\begin{equation}\label{SP_MTLL}
   \lim_{\eps\rightarrow 0}\eps \log_e
\tau_{c}=\lim_{\rho \rightarrow 0}\rho \log_e \tau_{c}=2.
\end{equation}

\begin{figure}
\centering \resizebox{!}{9cm}{\includegraphics{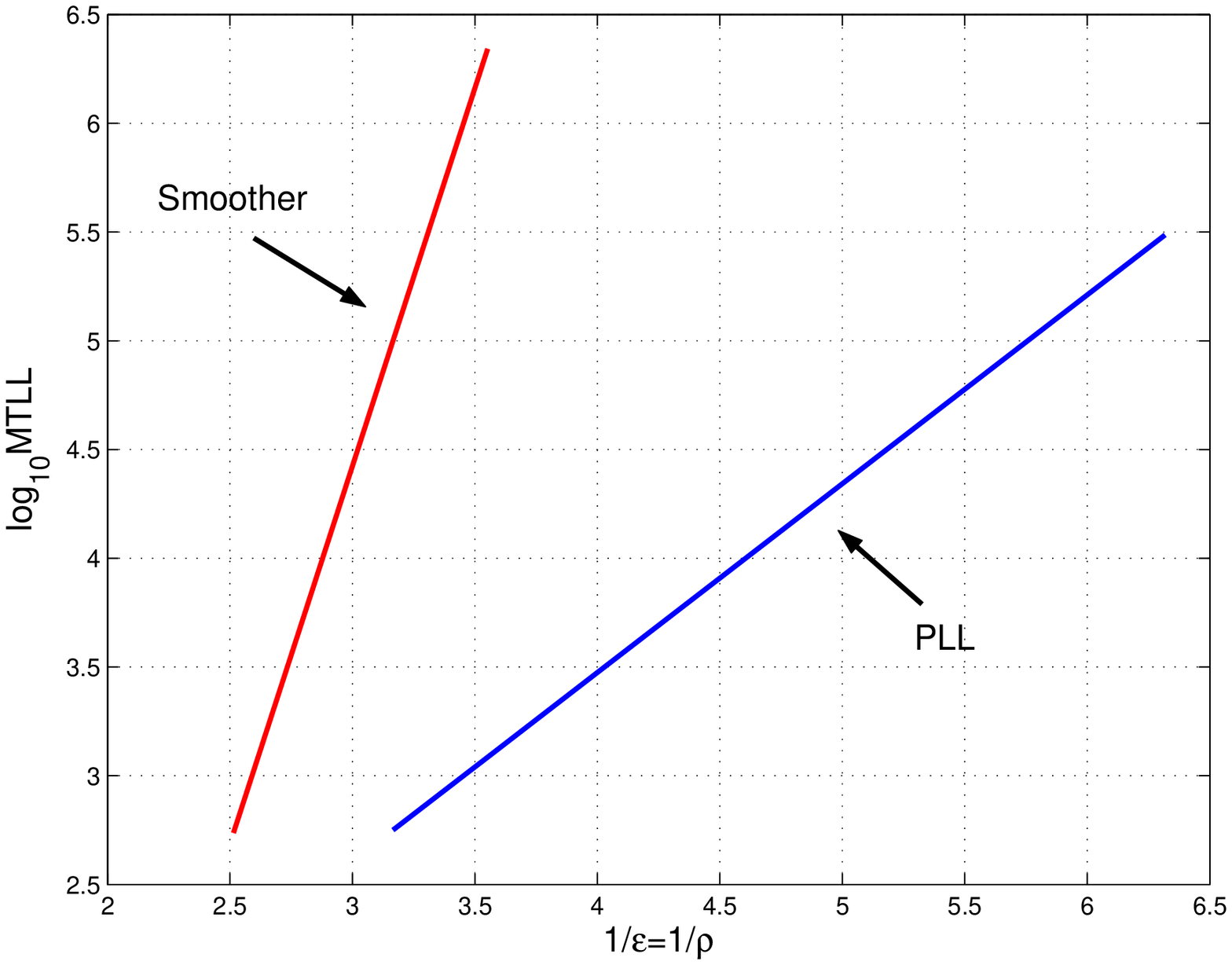}}
\caption{\label{MTLL1}The MTLL in the first order optimal MNE
estimator and the causal PLL as a function of $1/\eps$.  }
\end{figure}

\begin{figure}
\centering \resizebox{!}{9cm}{\includegraphics{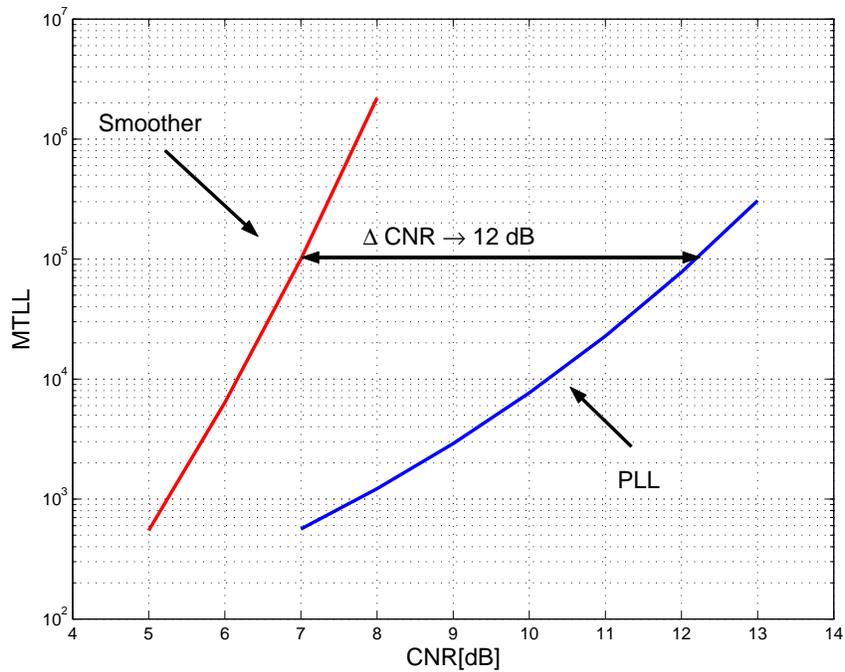}}
\caption{\label{MTLL_1_Ord_exp}The MTLL in the first order optimal
MNE estimator and the causal PLL as a function of the CNR.  }
\end{figure}

Note that in first order estimators the MTLL in the non-causal
estimator \eqref{ourMTLL} is significantly longer than that in the
causal estimator \eqref{SP_MTLL}. This implies that the CNR values
in the smoother are smaller than in the filter, but give a MTLL
identical to that in the filter. Denoting by $\eps_{c},\ \eps_{nc}$
the values of $\eps$ in the filter and smoother, respectively, and
requiring identical MTLLs,  lead to
\begin{eqnarray}
\tau_{c}=\tau_{nc} \Rightarrow
\frac{2}{\eps_{c}}=\frac{8}{\eps_{nc}}.
\end{eqnarray}
Denoting by $\mbox{CNR}_c[\mbox{dB}], \ \mbox{CNR}_{nc}[\mbox{dB}]$
the CNR in the filter and smoother, respectively, and using the
simple relation between $\eps$ and the CNR, lead to
\begin{eqnarray}\label{GapIndB}
\mbox{CNR}_{c}[\mbox{dB}]-\mbox{CNR}_{nc}[\mbox{dB}]=10 \cdot 2
\log_{10}(8/2)\approx 12 \mbox{dB}.
\end{eqnarray}
Thus, there exists a 12dB performance gap in the MTLL between the
estimators in terms of CNR. The MTLL in first order non-causal MNE
smoother and causal EKF are given in Figures \ref{MTLL1} and
\ref{MTLL_1_Ord_exp}. The pre-exponential term in the plots is
arbitrary.

Next we  consider the more complex, and more realistic, case of a
second order phase model \cite{Snyder}, in which the phase is
modelled as an integral over a Brownian motion. The signal model is
\begin{equation}\label{SecondDriving}
   \dot{\x}=\vm{A}'\x+\vm{B}'\dot{\w},
\end{equation}
where
\begin{eqnarray*}
\vm{A}'=\left(\begin{array}{cc}
      0 & 1 \\
      0 & 0
    \end{array}\right),\quad\vm{B}'= \left(\begin{array}{cc}
      0 & 0 \\
      0 & 1
    \end{array}\right).
\end{eqnarray*}
The measurements model is
\begin{equation}\label{sys2nd}
    \y=\vm{h}(\x) + \rho\, \dot{\vm{v}},
\end{equation}
where $\x(t)$ and $\w(t)$ take values in $ \mathbb{R}^2$. The
system \eqref{SecondDriving}, \eqref{sys2nd} is scaled to the form
of \eqref{GenralScaledSystemx},\eqref{GenralScaledSystemy} with
$\vm{A}=\vm{A}',\ \vm{B}=\vm{B}'$, and $\eps=\rho^{3/2}$. Note
that in this case the CNR equals $\rho^{-2}/2$, as in the case of
the first order system.

In the second order case, the asymptotic expression for the MTLL in
the smoother \eqref{generalMTLL} becomes
\begin{equation}\label{MTLLSecond}
\lim_{\eps \rightarrow 0}\eps \log_e \tau_{nc}= \mathop{\inf}
\limits_{e(\cdot) \in \mathscr{C}^2(t_0)} \frac{1}{2}\int_0^{T}\left
[ 4\sin^2 \left( \frac{e}{2} \right)+\ddot{e}^2 \right ] \,dt,
\end{equation}

variational problem on the right side of \eqref{generalMTLL} as no
analytical solution. An approximate solution is obtained numerically
about the characteristics of the Hamilton-Jacobi-Belman equation
\cite{Schuss_book,Benzi82}. This leads to the asymptotic limit of
the MTLL in the second order non-causal estimator as
\begin{equation}\label{ourMTLL2}
\lim_{\eps\rightarrow 0}\eps \log_e \tau_{nc}= \lim_{\rho
\rightarrow 0}\rho^{3/2} \log_e \tau_{nc}=5.
\end{equation}

\begin{figure}
\centering \resizebox{!}{9cm}{\includegraphics{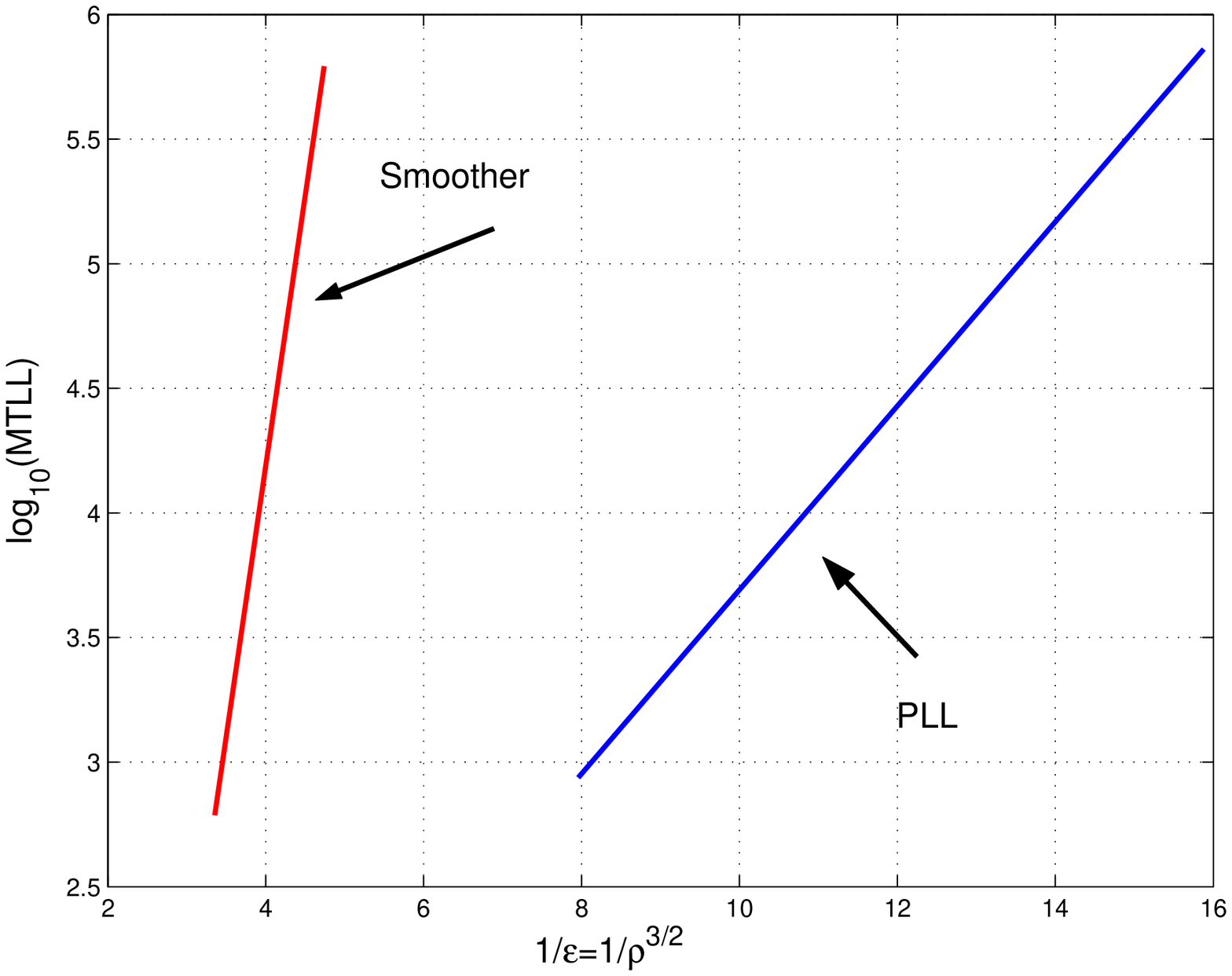}}
\caption{\label{MTLL2}The MTLL in the second order optimal MNE
estimator and the causal PLL as a function of $1/\eps$.}
\end{figure}

\begin{figure}
\centering \resizebox{!}{9cm}{\includegraphics{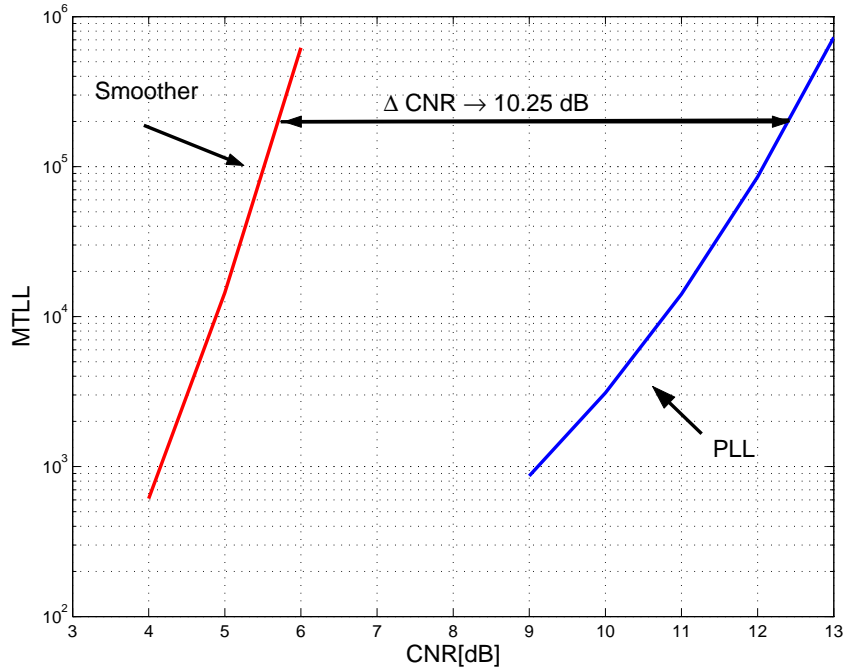}}
\caption{\label{MTLL_2_Ord_exp}The MTLL in the second order optimal
MNE estimator and the causal PLL as a function of CNR.}
\end{figure}

We note that the  EKF corresponding to the second order system
\eqref{SecondDriving}, \eqref{sys2nd} has the error equations
\cite{Snyder},\cite{Benzi82}
\begin{eqnarray}\label{error2nd}
\begin{array}{lll}
&&\\
 \dot{e}&=&\ds\frac{1}{2}\varphi-\sin e -\sqrt{2\eps'}\,\dot{v}\\
&&\\
 \dot{\varphi}&=& - \sin e - \sqrt{2\eps'}\,\dot{v} +
\sqrt{2\eps'}\,\dot{w}, \\
&&
\end{array}
\end{eqnarray}
where $\eps'=\eps/\sqrt{2}$. The "eikonal" equation \cite{Benzi82}
for the quasi-potential $\Phi$ corresponding to \eqref{error2nd} is
\begin{equation}\label{eiconal}
   \left(\frac{1}{2}\varphi-\sin e \right)\Phi_e -\sin e\,
   \Phi_\varphi+\Phi_{e}^2+2\Phi_{e}\Phi_{\varphi}+2\Phi_{\varphi}^2=0.
\end{equation}
The function $\Phi$ is evaluated numerically on the characteristics
of \eqref{eiconal} \cite{Benzi82}. This procedure leads to the value
of the quasi-potential at the unstable equilibrium point
$\Phi(\varphi=0,e=\pi)=0.6 $, so the asymptotic MTLL in the second
order causal estimator is
\begin{equation}\label{MTLL2ndEKF}
\lim_{\eps\rightarrow 0}\eps \log_e \tau_c= \lim_{\rho \rightarrow
0}\rho^{3/2} \log_e \tau_c=0.85.
\end{equation}
Note that similarly to the first order case, the MTLL in the second
order smoother \eqref{ourMTLL2} is significantly longer than that in
the causal PLL \eqref{MTLL2ndEKF}. The CNR gap in this case is
\begin{eqnarray}\label{GapIndB}
\mbox{CNR}_{c}[\mbox{dB}]-\mbox{CNR}_{nc}[\mbox{dB}]=40/3
\log_{10}(5/0.85)\approx 10.25 \mbox{dB}.
\end{eqnarray}
The MTLL in the non-causal estimator and the causal EKF in the
second order system are given in Figures \ref{MTLL2} and
\ref{MTLL_2_Ord_exp}. The pre-exponential term in the plots is
arbitrary.

\section{Discussion and conclusions}
The significant advantage of the optimal smoother over the causal
PLL defies intuition. It was customary to think that the MTLL
advantage of the smoother is linked and proportional to the MSE
advantage of the smoother \cite{VanTrees}. We argue that the MSE in
an estimator is not related to the MTLL. In order to demonstrate
this idea we introduce the following error equation
\begin{eqnarray}\label{CounterExample1}
\dot{e}=-2\sin \frac{e}{2} + \sqep \, \dot{\tilde{v}}.
\end{eqnarray}
The MSE in the linearized version of \eqref{CounterExample1} is
identical to that in linearized version of the system
\begin{eqnarray}\label{CounterExample2}
\dot{e}=-\sin e + \sqep \, \dot{\tilde{v}}.
\end{eqnarray}
However, due to the difference in the potential barrier in the above
systems, the asymptotical MTLL, $\tau_1$, in the system
\eqref{CounterExample1} satisfies
\begin{eqnarray}
\lim_{\eps \to 0}\eps \log_e \tau_1=\frac{1}{2}\mathop{\inf}
\limits_{\{e(0)=0, \ e(T')=2\pi\}} \int_0^{T'} \left ( \dot{e}+
2\sin \frac{e}{2} \right )^2\, dt=8,
\end{eqnarray}
while the MTLL, $\tau_2$, in the system \eqref{CounterExample2}
satisfies
\begin{eqnarray}
\lim_{\eps \to 0}\eps \log_e \tau_2=4,
\end{eqnarray}
and is significantly shorter than $\tau_1$. Thus, we conclude that
the MTLL advantage of the optimal smoother over the PLL is due to an
entirely different  functional and cannot be predicted by the MSE
advantage of the smoother over the filter.

There is a fundamental mathematical difference between the causal
and the non-causal cases. Both problems involve the minimization of
a functional, similar to that of the Wentzell-Freidlin theory for
causal systems. There is, however, a difference between the
functionals in the two theories. While the  functional in the causal
case vanishes along the exiting trajectories from the boundary of
the domain of attraction of the locked state to the next locked
state \cite{WentzellFreidlin}, in the non-causal case the functional
vanishes only at the locked states, so it has to be computed along
the entire slip trajectory.

\end{document}